\documentclass[12pt]{article}

\usepackage{amsmath}
\usepackage{amsfonts}
\usepackage{amssymb}
\newtheorem{theorem}{Theorem}[section]
\newtheorem{lemma}{Lemma}[section]
\newtheorem{corollary}{Corollary}[section]

\newtheorem{remark}{Remark}[section]

\begin{document}
\title{   Hardy Space Decompositions of $L^p(\mathbb{R}^n)$ for $0<p<1$ with Rational Approximation  }
\author{  Guan-Tie Deng, Hai-Chou Li,  Tao Qian}

\author{
Guan-Tie~Deng \thanks
{School of  Mathematical Sciences, Beijing
Normal University, Beijing, 100875, China.\ Email:denggt@bnu.edu.cn. This
work was partially supported by NSFC (Grant  11271045) and by SRFDP
(Grant 20100003110004)},\
Hai-Chou~Li \thanks{Corresponding author. College of Mathematics and informatics, South China Agricultural University, Guangzhou, China. Email: hcl2016@scau.edu.cn (or, lihaichou@126.com).},\
        Tao~Qian \thanks{Department of Mathematics, University of Macau, Macao (Via Hong
Kong). Email: fsttq@umac.mo. The work was partially supported by
Research Grant of University of Macau, MYRG115(Y1-L4)-FST13-QT, Macao Government FDCT 098/2012/A3}
}
\date{}

\maketitle
\begin{center}
\begin{minipage}{120mm}
\begin{center}{\bf Abstract}\end{center}
{This paper aims to obtain decompositions of higher dimensional $L^p(\mathbb{R}^n)$ functions into  sums of non-tangential boundary limits of the corresponding Hardy space functions on tubes for the index range $0<p<1$. In the one-dimensional case, Deng and Qian \cite{DQ} recently obtained such Hardy space decomposition result: for any function $f\in L^p(\mathbb{R}),\ 0<p<1$, there exist functions $f_1$ and $f_2$ such that $f=f_1+f_2$, where $f_1$ and $f_2$ are, respectively, the non-tangential boundary limits of some Hardy space functions in the upper-half and lower-half planes. In the present paper,  we generalize the one-dimensional Hardy space decomposition result to the higher dimensions, and discuss the uniqueness issue of such decomposition. }

{\bf Key words}:\ \ Hardy spaces on tube domain,  decomposition, rational approximation
\end{minipage}
\end{center}

 \section { Introduction}

  We  begin with a survey on  the case of one dimensional Hardy spaces $H^p(\mathbb{C}^\pm)$ for $0<p\leq\infty$. For the classical case $p=2$,  the famous Paley-Wiener Theorem states that, for an $L^2(\mathbb{R})$ function $f,$ it is the non-tangential boundary limit of a function in $H^2(\mathbb{C}^+)$ if and only if ${\rm supp}\hat{f}\subset [0,\infty),$ the latter being equivalent with
\begin{eqnarray}\label{spectrum characterization for conjugate H2}\hat{f}=\chi_{[0,\infty)}\hat{f},\end{eqnarray}
where $\chi_{[0,\infty)}$ is the characteristic function of $[0,\infty)$.
Similarly,
$f\in H^2(\mathbb{C}^-)$  is equivalent with
\begin{eqnarray}\label{spectrum characterization for H2}\hat{f}=\chi_{(-\infty,0]}\hat{f}.\end{eqnarray}

From  the relations (\ref{spectrum characterization for conjugate H2}) and (\ref{spectrum characterization for H2}), a canonical decomposition
result is obtained as
 $$f=f^++f^-,$$
  for $f\in L^2(\mathbb{R}),$ where $f^+=(\chi_{[0,\infty)}\hat{f})^\vee$ and $f^-=(\chi_{(-\infty, 0]}\hat{f})^\vee$
belong to, respectively, the Hardy spaces $H^2(\mathbb{C}^+)$ and $H^2(\mathbb{C}^-).$  In the present paper we  call such decomposition as Hardy space decomposition.  It can  be alternatively obtained through the Hilbert transformation in view of the Plemelj formula in relation to the Cauchy integral transform. For the relevant knowledge we refer the reader to, for instance, \cite{Gar1} and \cite{Q2}.

The significance of the $L^2(\mathbb{R})$ space decomposition into the Hardy spaces $H^2(\mathbb{C}^\pm)$ lays on the fact that
the $H^2$ functions have much better properties than the $L^2$-functions. The functions $f^\pm$ are non-tangential boundary
limits of $H^2(\mathbb{C}^\pm)$ functions, the latter being analytically defined in their respective domains. We note that a function in
$L^p$ is only a.e. determined, that is, if $f=g$, a.e., then $f$ and $g$ are considered as  the same function in $L^p.$ In this regard one cannot assume a general $L^p$-function to have any  smoothness. On the other hand, any two a.e. identical $L^p$
functions correspond to the same Hardy space decomposition, the latter, being unique, having all kinds of smoothness in their domains of definition. The Hardy $H^2$ spaces are reproducing kernel Hilbert spaces. In particular, Cauchy's theory and techniques are available to Hardy space functions. Furthermore, there is an isometric correspondence between the Hardy space functions and their non-tangential boundary limits. Based on what just mentioned,
 analysis on the $L^2$-functions can be reduced to that of  their Hardy space components. Such treatment of
$L^2$-functions has found significant applications in both the pure and applied mathematics.  In particular, some
demonstrative results in signal analysis have recently been obtained ( \cite{DQY}, \cite{PA1},
\cite{PMA},).

There have been studies on Hardy space decomposition  for extended index range of $p$. Both the above mentioned Paley-Wiener Theorem related  Fourier spectrum characterization and the Plemelj formula  can be extended to Hardy
$H^p(\mathbb{R})$ spaces with  $p\ne 2.$ Systematic studies on the spectrum properties as well as the $L^p$
decomposition are carried out in the works of Qian et al. \cite{Q2} among others. Their work is summarized as  $f\in H^p(\mathbb{C}^+)$ if and only if $f\in L^p(\mathbb{R})$ and in the distributional sense ${\rm supp}\hat{f}\subset [0,\infty)$ for all $1\leq p\leq \infty.$ Comparatively the Plemelj formula
approach is more applicable  than Fourier transformation, for on the $L^p(\mathbb{R}), \ p>2,$ spaces Fourier transforms are distributions. As a consequence, the Hardy space decomposition results also hold for $f\in L^p(\mathbb{R}), \ 1<p<\infty.$ There also exists an analogous theory on the unit circle for $1< p\leq \infty.$ The latter context corresponds to  Fourier series other than Fourier transforms (\cite{Gar1}, \cite{CJA}).

Turning to  the case of $0<p\leq 1,$ due to loss of integrability neither on the real line, nor on the unit circle, the Plemelj formula, or the Hilbert transformation are available. The references \cite{CJA} and \cite{DQ} study the Hardy space decomposition of $L^p(\mathbb{R})$ for $0<p<1.$ There is no Fourier transformation theory for functions in such $L^p(\mathbb{R})$  as they are even not distributions. One, however, can still have the Hardy space decompositions. The reference \cite{CJA} uses the real analysis methods of harmonic analysis, making use of a dense subclass of the $L^p(\mathbb{R})$-functions with vanishing moment conditions and the Hilbert transforms. In contrast,  \cite{DQ} uses the complex analysis methods, and, in particular, rational functions approximation to achieve Hardy space decompositions. The methods in \cite{DQ} are direct and constructive. The following Theorem A is obtained by Deng and Qian in \cite{DQ}.  \\
 {\bf  Theorem A }\ (\cite{DQ}) \ Suppose that   $0<p<1$ and  $ f\in L^p(\mathbb{R})$.
  Then,  there exist a positive constant $A_p$ and  two  sequences of rational functions $\{P_k\}$ and $\{Q_k\}$ such that
  $ P_k\in H^p(\mathbb{C}_+)$, $ Q_k\in H^p(\mathbb{C}_-)$
  and   $$
  \sum_{k=1}^{\infty} \left (\|P_k\|^p_{H^p_{+}}+\|Q_k\|^p_{H^p_{-}}\right )\leq A_p\|f\|_p^p ,
  $$
  $$
     \lim_{n\rightarrow \infty}||f-\sum_{k=1}^{n}(P_k+Q_k)||_p
    =0.
  $$
  Moreover,
  $$
  g(z)=\sum_{k=1}^{\infty}P_k(z)\in H^p(\mathbb{C}_+) ,\ \ h(z)=\sum_{k=1}^{\infty}Q_k(z)\in H^p(\mathbb{C}_-),
  $$
  and
   $g(x)$ and $ h(x)$ are the non-tangential boundary
 limits   of  functions for $g\in H^p(\mathbb{C}_+)$ and $h\in
  H^p(\mathbb{C}_-)$, respectively, $f(x)=g(x)+h(x)$ almost
  everywhere, and
  $$\|f\|_p^p\leq  \|g\|_p^p+\|h\|_p^p\leq A_p\|f\|_p^p . $$
    That is, in the sense of $L^p(\mathbb{R})$,
  $$
  L^p(\mathbb{R})=H^{p}_+(\mathbb{R})+H^p_-(\mathbb{R}),
  $$
 where $H_+(\mathbb{R})$ and $H_-(\mathbb{R})$  denote the set of non-tangential boundary limits of functions in the Hardy
 spaces
 $H^p(\mathbb{C}_+)$ and $H^p(\mathbb{C}_-)$, respectively.

Recently we, in \cite{LDQ1}, obtain the  Hardy space decomposition of $L^p$ for $0<p\leq 1$ on the unite circle $\partial \mathbb{D}$ by using polynomial approximation other than general rational approximation. The result is stated as follows. Denoting  by  $ L^p_I(\partial \mathbb{D}) $ and  $ L^p_O(\partial \mathbb{D})$ the closed subspaces of $L^p(\partial \mathbb{D}),$ consisting of respectively, the non-tangential boundary limits of the functions of $H^p(\mathbb{D})$ and $H^p(\mathbb{D}_O),$   then we have
\[ L^p(\partial \mathbb{D})=L^p_I(\partial \mathbb{D}) + L^p_O(\partial \mathbb{D}),\]
where the right-hand-side is not a direct sum. As a matter of fact, the intersection $L^p_I(\partial \mathbb{D})$ and
$ L^p_O(\partial \mathbb{D})$  contains non zero functions. The work on the unit circle exposes the particular features adaptable to higher dimensions.

All the Hardy space  results mentioned above for dimension one can be generalized to  dimension n in the setting of Hardy spaces on tubes with correspondingly the right notions. In fact, \cite{SW} already contain some basic results, mostly for $p=2$, while \cite{LDQ} gives a systematic treatment for general indices $p\in[1,\infty]$, including Fourier spectrum characterization of Hardy spaces on tubes, the Cauchy integral and Poisson integral representation of the Hardy space functions, the Plemelj formulas in relation to Hilbert transforms, and the Hardy space decompositions of functions in the $L^p(\mathbb{R}^n)$ spaces. The purpose of this article is to prove the Hardy space decomposition of the $L^p(\mathbb{R}^n)$ space functions for $p\in (0,1)$. In doing so neither the Cauchy integral formula nor the Fourier transformation can be directly used, for the functions defined on $\mathbb{R}^n$ are lack of integrability. They are even not distributions.

In the present paper, inspired by the idea of \cite{DQ}, with the rational approximation method, we obtain decompositions of functions in $L^p(\mathbb{R})$ for $0<p<1$ into sums of boundary limits of the corresponding Hardy space functions on tubes,  $H^p(T_{\Gamma_{\sigma_j}})\ (j=1,2,...,2^n)$, through the rational functions approximation. The idea of using rational approximation is motivated by the studies \cite{WJL} of Takenaka-Malmquist systems.

We will discuss the non-uniqueness of the Hardy space decomposition  via rational approximation method.  We conclude that the sum is not a direct sum. In fact, there is a non-trivial intersection of all those summed spaces. The intersection of those spaces $H_{\sigma_j}^{p}(\mathbb{R}^n)\ ( j=1,2,...,2^n)$ is identical with the $L^p$ closure of the set of functions $f\in L^p(\mathbb{R}^n)$ of  the following form
$$f(x_1,x_2,...,x_n)= \frac{P(x_1,...,x_n)}{\prod\limits_{k=1}^n\prod\limits_{j=1}^m(x_{k}-a_{kj})},$$
where $P(x_1,...,x_n)= \sum\limits_{s_1=0}^{l_1}\cdot\cdot\cdot\sum\limits_{s_n=0}^{l_n}a_{s_1...s_n}x_1^{s_1}\cdot\cdot\cdot
x_n^{s_n},\ l= \max\{l_1,...,l_n\}, $ $(m-l)p>1,$ and  $a_{kj}\neq a_{km} \in \mathbb{R}$ as $ j\neq m$ for $ k=1,2,...,n.$

 The writing plan of this paper is as follows:  In \S\ 2, some basic definitions and notations are given. In \S\ 3 we devote to establishing the higher dimensional Hardy space decomposition of $L^p(\mathbb{R}^{n}), \ 0< p< 1.$ The decomposition is   a sum of  boundary limit functions of  Hardy spaces on tubes, $H^{p}(T_{\Gamma_{\sigma_k}})$, for all $ k = 1,2,..., 2^{n}.$ In
 \S\ 4, we discuss the uniqueness of such Hardy space decomposition.

 \section {Preliminary Knowledge}

 In this section, we introduce some useful basic definitions and notions. For more information, see e.g.  \cite{Gar1} and \cite{SW}.
The classical Hardy spaces $ H^p(\mathbb{C}_k) ,\ 0<p< +\infty,\ k=\pm 1 $,
consists of the functions $f$ analytic in the
half plane $ \mathbb{C}_k=\{z=x+iy: k y>0\}.$ They are Banach spaces for $1\leq p<\infty$ under the norms
$$
 \|f\|_{H^p_k}=\sup_{k y>0}\left
(\int_{-\infty}^{\infty}|f(x+iy)|^pdx\right )^{\frac{1}{p}};
    $$
    and complete metric spaces for $0<p<1$ under the metric functions
    $$ d(f,g)=\sup_{k y>0}\int_{-\infty}^{\infty}|f(x+iy)|^pdx.$$
Let $B$ be an open subset of $\mathbb{R}^{n}$. Then the tube $T_B$ with base $B\subset \mathbb{R}^n$ is the set
$$T_B=\{z=x+iy\in \mathbb{C}^n:\ x\in \mathbb{R}^{n},\ y\in B\}.$$
For example, when $n=1$, the classical upper-half complex plane $\mathbb{C}^{+}$ and lower-half complex plane $\mathbb{C}^{-}$ are the tubes in $\mathbb{C}$ with the base $B_+=\{y\in\mathbb{R}:\ y>0\}$ and the base $B_-=\{y\in\mathbb{R}:\ y<0\}$, respectively. That is,
$\mathbb{C}^{+}=T_{B_+}=\{z=x+iy:\ x\in \mathbb{R},\ y>0\}$ and $\mathbb{C}^{-}=T_{B_-}=\{z=x+iy:\ x\in \mathbb{R},\ y<0\}$.
 Obviously, the tube $T_B$ are generalizations of $\mathbb{C}^{+}$ and  $\mathbb{C}^{-}$.

It is known that $n$-dimensional real Euclidean space $\mathbb{R}^n$ has $2^{n}$ octants. To denote the octants we adopt the following notations.
Let $\sigma_k= (\sigma_k(1),\sigma_k(2),...,\sigma_k(n))$, where $ \sigma_k(j)=+1\ {\rm or} -1$ for $ k = 1,2,..., 2^{n}$. Then $2^{n}$ octants of $\mathbb{R}^n$ are denoted by $\Gamma_{\sigma_k},\ \ k = 1,2,..., 2^{n}$, where
$$\Gamma_{\sigma_k}=\{y=(y_1,y_2,...,y_n)\in\mathbb{R}^{n}: \ \sigma_k (j)y_j>0, \ j=1,2,...,n\}.$$
In this paper, we denote $\Gamma_{\sigma_1}$ the first octant of $\mathbb{R}^n$, that is
$$\Gamma_{\sigma_1}=\{y=(y_1,y_2,...,y_n)\in\mathbb{R}^{n}: \ y_j>0, \ j=1,2,...,n\}.$$
Correspondingly,
$\mathbb{C}^{n}$  can be decomposed into $2^{n}$
tubes, denoted by $T_{\Gamma_{\sigma_k}},\ k = 1,2,..., 2^{n}$. That is
$$T_{\Gamma_{\sigma_k}}=\{z=x+iy\in\mathbb{C}^{n}:\ x\in\mathbb{R}^{n},\ y\in\Gamma_{\sigma_k}\}.$$

A function $F(z)$ is said to belong to the space $H^{p}(T_{B}),\ 0< p<\infty,$ if it is holomorphic in the  tube $T_B$, and
satisfies
$$
\|F\|_{H^p}=\sup \left \{\left ( \int_{\mathbb{R}^{n}}|F(x+iy)|^{p}dx\right )^{\frac{1}{p}} :\ y\in B \right \} < \infty.
$$
Hence,
$$H^{p}(T_{B})=\{F: F \mbox{ holomorphic on  }\  T_{B} \mbox{ and }\   \|F\|_{H^p} < \infty \}.$$
The spaces $H^{p}(T_{\Gamma_{\sigma_k}})$ are defined through replacing $B$ by $\Gamma_{\sigma_k},\ k=1,...,2^n.$

 A function, $f$, defined in tube $T_\Gamma$, is said to have non-tangential boundary limit (NTBL) $l$ in each component of the
 variable  at $x_0\in \mathbb{R}^{n}$ if $f(z)=f(x+iy)=f(x_1+iy_1,...,x_n+iy_n)$ tends to $l$ as the point $z=(x_1,y_1;
 x_2,y_2;...; x_n,y_n)$ tends to $x_0=(x_0^{(1)},0; x_0^{(2)},0;...; x_0^{(n)},0)$ within the Cartesian product
$$\gamma _\alpha(x_0)= \Gamma_{\alpha_1}(x_0^{(1)})\times \Gamma_{\alpha_2}(x_0^{(2)})\times \cdot\cdot\cdot\times
\Gamma_{\alpha_n}(x_0^{(n)})\subset T_\Gamma,$$
for all $n$-tuples $\alpha= (\alpha_1,\alpha_2,...,\alpha_n)$ of positive real numbers, where
$$\Gamma_{\alpha_j}(x_0^{(j)})=\left\{(x_j,y_j)\in \mathbb{C}^{+}: \ |x_j-x_0^{(j)}|< \alpha_jy_j\right\} ,\ j=1,2,...,n.$$
As an important property of the Hardy spaces, it is shown that if $f$ is a function in a Hardy space $H^p, 0<p< \infty,$ then
for almost all $x_0, \ f$ has NTBL (\cite{SW}).

Since the mapping that maps the functions in the Hardy spaces to their NTBLs is an isometric isomorphism, we denote by
$H_{\sigma_k}^{p}(\mathbb{R}^{n})$ for $0<p<1$ the NTBLs of the functions in  $H^{p}(T_{\Gamma_{\sigma_k}})$, that is
$$H_{\sigma_k}^{p}(\mathbb{R}^{n})=\{f: \ f \mbox{ is the NTBL of a function in } H^{p}(T_{\Gamma_{\sigma_k}})\}$$
for all  $ k = 1,2,..., 2^{n}$. The non-tangential
boundary limit of $F(z)\in H^{p}(T_{\Gamma_{\sigma_k}})$ as $y \rightarrow 0$ in the tube are denoted by
\begin{equation}\label{57}
F_{\sigma_k}(x)= \lim \limits _{y \in \Gamma_{\sigma_k}, y\rightarrow 0 }F(x+iy)
= \lim \limits _{ \sigma_k(1)y_1\rightarrow 0^+,..., \sigma_k(n)y_n\rightarrow 0^+ }F(x_1+iy_1,...,x_n+iy_n).
\end{equation}

 \section{ Hardy Space Decomposition of $L^p(\mathbb{R}^n)$}

As previously mentioned in the introduction, for one dimension, in \cite{DQ}, the authors use the rational approximation method to obtain the Hardy space decomposition for the range $0<p<1$. As is well known, rational approximation has a long history, and is naturally related to complex approximation \cite{WJL}. As obtained in \cite{DQ} through a rational approximation method, for a given real-valued function $f\in L^p(\mathbb{R}),\ 0<p<1$,  there exists the relation $f=f_1+f_2$, where $f_+$ and $f_-$ are the non-tangential boundary limit functions of some analytic $H^p$ functions in, respectively, the upper-half and the lower-half complex planes \cite{DQ}.
 Precisely, the analytic Hardy space function  $f_+\in H^p_+$,  and thus its boundary limit function as well, are defined through a sequence of rational functions whose poles are in the lower-half plane, and   $f_-\in H^p_-$, through a sequence of rational functions whose poles are in the upper-half plane. We note that the Hardy spaces decompositions for functions in $L^p(\mathbb{R}), \ 0<p<1,$ are not unique. This amounts to saying that the intersection  $H^p_+ \bigcap H^p_-$ a non-empty set.

In this section by using a higher dimensional rational approximation method, we will generalize the above type of Hardy space decomposition  of $L^{p}(\mathbb{R})$ to higher dimensional $L^{p}(\mathbb{R}^{n}),\ 0<p<1,$ and obtain Theorem 3.1. Specifically, for any real-valued function $f\in L^p(\mathbb{R}^n),\ 0<p<1$, it is proved to  have the Hardy space decomposition  $f(x)=\sum\limits_{j=1}^{2^n}f_{\sigma_j}(x),$ where for each $j$, $f_{\sigma_j}(x) \ (j=1,2,...,2^n)$ is the non-tangential boundary limit of some $H^p(T_{\Gamma_{\sigma_j}})$-function.  In fact, each analytic Hardy space function $f_{\sigma_j}(x)\in H^p_{\sigma_j}$, and thus its boundary limit function as well, may be approximated by a sequence of rational $L^{p}(\mathbb{R}^{n})$-functions whose poles are not in the octant $T_{\Gamma_{\sigma_j}}$. We will call such rational functions rational atoms.

  \begin{theorem}\ \  Suppose that $f\in L^{p}(\mathbb{R}^{n}),\ 0<p< 1$. Then there exist $2^{n}$ sequences of rational
  functions $\{R_{{k\sigma}_{j}}(z)\}\in H^{p}(T_{\Gamma_{\sigma_j}}),$ and $ f_{\sigma_j}(z)\in H^{p}(T_{\Gamma_{\sigma_j}}),
  \ j=1,2,...,2^{n},$ such that the following properties hold\\
(i)
 \begin{equation}\label{2}
   \sum\limits _{k=1}^{\infty}\sum\limits_{j=1}^{2^n}||R_{{k\sigma}_{j}}||^p_{H^p_{\sigma_j}}\leq A_{np}\|f\|_p^p,
  \end{equation}
  where $A_{np}$ is a constant only depending on $(n,p)$;\\
(ii)
 \begin{equation}\label{3}
    \lim_{m\rightarrow \infty}||f-\sum_{k=1}^{m}\sum\limits_{j=1}^{2^n}R_{{k\sigma}_{j}}||_p
  =0;
  \end{equation}
(iii)
 \begin{equation}\label{4}
  \lim_{m\rightarrow \infty}||f_{\sigma_j}-\sum_{k=1}^{m}R_{{k\sigma}_{j}}||_{H^p}=0 ,
  \end{equation}
for all $j=1,2,...,2^{n};$ \\
(iv) $f_{\sigma_j}(x)$ are the non-tangential boundary limits of functions $f_{\sigma_j}(z)$, and
$$f(x)=\sum\limits_{j=1}^{2^n}f_{\sigma_j}(x), \ \ \ \ \  a.e.\ \  x\in \mathbb{R}^n;$$
(v)
$$\|f\|^{p}_{p} \leq \sum\limits_{j=1}^{2^n}\|f_{\sigma_j}\|^{p}_{p}\leq  A_{np}\|f\|_p^p.$$
(vi) In summary,
$$L^p(\mathbb{R}^{n})=\sum\limits_{j=1}^{2^n}H_{\sigma_j}^{p}(\mathbb{R}^n).$$
 \end{theorem}

 \begin{remark}\ \ For the non-uniqueness issue of the Hardy space decomposition,we will show that $\bigcap\limits_{j=1}^{2^n} H_{\sigma_j}^p(\mathbb{R}^n)$ is a non-empty set.  We will prove it  in \S\ 4.
\end{remark}

In order to prove  Theorem 3.1, we need the  following lemmas. We note that the proof of Theorem 3.1 is at the end of this section.
\begin{lemma}
\ Let  $\Gamma_{\sigma_1}$ be the first octant of $\mathbb{R}^n$. Suppose that   $0<p<1$ and  $R$ is a rational
function with the  form
$$R(z)=\frac{P(z)}{Q(z)},$$
where $P(z)$ is a polynomial of $z$, and $ Q(z)=Q_1(z_1)Q_2(z_2)\cdot\cdot\cdot\cdot Q_n(z_n)$, $Q(z_k)$ is a polynomial of $z_k,\
k=1,2,...,n$.
 If   $ R\in L^p(\mathbb{R}^n)$ and $R(z)$ is holomorphic in $T_{\Gamma_{\sigma_1}}$, then
$ R(z)\in H^p(T_{\Gamma_{\sigma_1}}).$
\end{lemma}
 Deng-Qian  \cite{DQ} proved the special case $n=1$ of Lemma 3.1, and obtained the following Lemma 3.1-1, which is needed for the proof of Lemma 3.1.

{\bf  Lemma  3.1-1 }(\cite{DQ}) \ Suppose that   $0<p<1$ and  $R$ is a rational
function with $ R\in L^p(\mathbb{R}).$  For $k=\pm 1$, if  $R(z)$ is analytic in the half plane $ \mathbb{C}_{k}$, then
$ R\in H^p(\mathbb{C}_k)$.

In order to prove Lemma 3.1,  we  need the following Lemma 3.1-2 which  obtained in our recently work \cite{LDQ2}:

{\bf  Lemma  3.1-2 }(\cite{LDQ2})\  If $f(z)\in H^p(T_{\Gamma_{\sigma_j}}),\ 0<p<\infty, \ j=1,2,...,2^n$, and $f(x)$ is the boundary limit of $f(z)$. Then
$\varphi(y)$ is continuous convex and bounded in $\overline{\Gamma}_{\sigma_j}$, moreover,
$$\|f\|^p_{H^p_{\sigma_j}}=\sup\limits _{y\in \Gamma_{\sigma_j}}\varphi(y)=\varphi(0,...,0)= \|f\|^p_{p},$$
where $\varphi(y)= \int_{\mathbb{R}^{n}}|f(x+iy)|^p dx,\ y\in \overline{\Gamma}_{\sigma_j},\ j=1,2,...,2^n.$\\
{\bf  Proof of Lemma 3.1 }\ \ We are to prove Lemma 3.1 by mathematical induction.

 When $n=1$,  Lemma 3.1  is just  Lemma 3.1-1.

Next, when $n>1,$ we assume that Lemma 3.1 holds for $n-1$. Take $ n_1=n-1,$ and fix $ t_n\in \mathbb{R}\backslash Z_0(Q_n)$,
where $Z_0(Q_n)=\{z_n:\ Q_n(z_n)=0\}$. We consider the function $r$ of $n-1$ real variables defined by
$$r(x_1,x_2...,x_{n-1})= R(x_1,...,x_{n-1},t_{n}),$$
where $(x_1,x_2...,x_{n-1})\in \mathbb{R}^{n-1}.$
The Fubini Theorem ensures that $r(x_1,x_2...,x_{n-1})$ belongs to $L^p(\mathbb{R}^{n-1})$
for almost all $ t_n\in \mathbb{R}\backslash Z_0(Q_n).$ Moreover, it is easy to see that the rational function
$r(z_1,...,z_{n-1})$ satisfies the assumptions of Lemma 3.1.

Therefore, by the induction hypothesis, we obtain that $r(z_1,...,z_{n-1})\in H^p(T_{\Gamma^{n-1}_{\sigma_1}}),$ where
$T_{\Gamma^{n-1}_{\sigma_1}}$ denotes the tube with the first octant of $\mathbb{R}^{n-1}$ as base.

By Lemma 3.1-2,
 $$\int_{\mathbb{R}^{n-1}}|r(x_1+iy_1,...,x_{n-1}+iy_{n-1})|^{p}\ dx_1\cdot\cdot\cdot dx_{n-1} \leq
 \int_{\mathbb{R}^{n-1}}|r(x_1,...,x_{n-1})|^{p}\ dx_1\cdot\cdot\cdot dx_{n-1}.$$
 That is
 \begin{equation}\label{52}
 \begin{split}
&\int_{\mathbb{R}^{n-1}}|R(x_1+iy_1,...,x_{n-1}+iy_{n-1},x_n)|^{p}\ dx_1\cdot\cdot\cdot dx_{n-1}
\\
&\leq
  \int_{\mathbb{R}^{n-1}}|R(x_1,...,x_{n-1},x_n)|^{p}\ dx_1\cdot\cdot\cdot dx_{n-1}.
\end{split}
\end{equation}
 Integrating both sides of the last inequality with respect to $x_n$, we have
 \begin{equation}\label{50}
\int_{\mathbb{R}^{n}}|R(x_1+iy_1,...,x_{n-1}+iy_{n-1}, x_n)|^{p}\ dx_1\cdot\cdot\cdot dx_{n} \leq
\int_{\mathbb{R}^{n}}|R(x_1,...,x_{n})|^{p}\ dx_1\cdot\cdot\cdot dx_{n}.
\end{equation}
 By Fubini Theorem,
  \begin{equation}\label{51}
\int_{\mathbb{R}^{n-1}}\left(\int_{\mathbb{R}}|R(x_1+iy_1,...,x_{n-1}+iy_{n-1}, x_n)|^{p}\ dx_{n}\right)dx_1\cdot\cdot\cdot
dx_{n-1} \leq \int_{\mathbb{R}^{n}}|R(x)|^{p}\ dx.
\end{equation}
 So, fix $(y_1,...,y_{n-1})$, for almost all $(x_1,...,x_{n-1})\in \mathbb{R}^{n-1}$, we have
 $$\int_{\mathbb{R}}|R(x_1+iy_1,...,x_{n-1}+iy_{n-1}, x_n)|^{p}\ dx_{n}< \infty.$$
 That is, rational function $R(x_1+iy_1,...,x_{n-1}+iy_{n-1}, x_n)$ as a function of $x_n$ belongs to $ L^p(\mathbb{R})$.
 Moreover, since
 $R(z_1,...,z_n)$ is holomorphic in $T_{\Gamma_{\sigma_1}}$,  $R(z_1,...,z_n)$ as a function of $z_n$ is also holomorphic in
 upper-half plane $\mathbb{C}^+$.
 Due to the result for $n=1$,  $R(z_1,...,z_{n-1}, z_n)$ as a function of $z_n$ is a member of $H^p(\mathbb{C}^+)$. By Lemma 3.1-2 again,
 $$\int_{\mathbb{R}}R(z_1,...,z_{n-1}, x_n+iy_n)|^{p}\ dx_{n} \leq \int_{\mathbb{R}}|R(z_1,...,z_{n-1},x_{n})|^{p}\ dx_{n}.$$
Integrating both sides of the last inequality with respect to $x_1,...,x_{n-1}$, we have
 \begin{equation}\label{52}
 \begin{split}
&\int_{\mathbb{R}^{n-1}}\left(\int_{\mathbb{R}}|R(x_1+iy_1,...,x_{n-1}+iy_{n-1}, x_n+iy_n)|^{p}\
dx_{n}\right)dx_1\cdot\cdot\cdot dx_{n-1}\\
&= \int_{\mathbb{R}^{n}}|R(x+iy)|^{p}\ dx\\
& \leq \int_{\mathbb{R}^{n-1}}\left(\int_{\mathbb{R}}|R(x_1+iy_1,...,x_{n-1}+iy_{n-1}, x_n)|^{p}\
dx_{n}\right)dx_1\cdot\cdot\cdot dx_{n-1}.
\end{split}
\end{equation}
Together with \eqref{51} and \eqref{52}, we obtain
$$\int_{\mathbb{R}^{n}}|R(x+iy)|^{p}\ dx \leq \int_{\mathbb{R}^{n}}|R(x)|^{p}\ dx.$$
Due to  the analyticity of the rational function $R(z)$ in $T_{\Gamma_{\sigma_1}}$, we get $R(z)\in
H^p(T_{\Gamma_{\sigma_1}}).$ The proofs for the other octants are similar. So the proof of Lemma 3.1 is complete.\\

Similarly to prove Lemma 3.1, we also can get the analogous results about the other octants as follows.

\begin{corollary}\ \
\ Let  $\Gamma_{\sigma_j}\ ( j=1,2,...,2^n)$ be  all the octants of $\mathbb{R}^n$. Suppose that   $0<p<1$ and  $R$ is a
rational function with the  form
$$R(z)=\frac{P(z)}{Q(z)},$$
where $P(z)$ is a polynomial of $z$, and $ Q(z)=Q_1(z_1)Q_2(z_2)\cdot\cdot\cdot\cdot Q_n(z_n)$, where for each $k$, $Q_k(z_k)$ is a polynomial of $z_k,\
k=1,2,...,n$.
 If   $ R\in L^p(\mathbb{R}^n)$ and $R(z)$ is holomorphic in $T_{\Gamma_{\sigma_j}}$, then
$ R(z)\in H^p(T_{\Gamma_{\sigma_j}}),$  $j=1,2,...,2^n.$
\end{corollary}

\begin{lemma}\ \  If  $0<p<1$,  $ f\in L^p(\mathbb{R}^{n})$, then, for $\varepsilon>0$,
 there exist a sequence of rational functions $\{R_k(x)\}$, $R_k\in \mathcal{A}$, such that
 \begin{equation}\label{6}
  \sum_{k=1}^{\infty }||R_k||^p_p\leq (1+\varepsilon )\|f\|_p^p,
 \end{equation}
and
 \begin{equation}\label{7}
   \lim_{n\rightarrow \infty}||f-\sum_{k=1}^{n}R_k||_p
  =0,
 \end{equation}
 where
 $$ \mathcal{A}= \left\{R(x)=\frac{P(x)}{(1+x_1^{2})^{l_1}\cdot\cdot\cdot (1+x_n^{2})^{l_n}}: \ 2l_j> {\rm deg}_jP,\
 j=1,2,...,n\right\},$$
 $P(x)=\sum_{(k)}\alpha_{k}  x_{1}^{k_1}
 \cdot\cdot\cdot x_{n}^{k_n}$ is a polynomial of $x=(x_1,x_2,...,x_n) \in \mathbb{R}^{n}$, $\alpha_{k}$ are constants,
 $k=(k_1,k_2,...,k_n)$, and  the notation $\sum_{(k)}$ means that
 $\sum_{(k)}=\sum\limits_{k_n=0}^{s_n}\cdot\cdot\cdot\sum\limits_{k_1=0}^{s_1},$ $ {\rm deg}_jP=s_j,\ j=1,2,...,n.$
\end{lemma}

{\bf Proof }\ \ \ We assume that $\|f\|_p>0,$ and let
 $$C_0(\mathbb{R}^{n})=\left\{f:\ f \ {\rm  is\ continuous \ in } \ \mathbb{R}^{n}, \
 \lim\limits_{|x|\rightarrow\infty}f(x)=0\right\}.$$
 It is obviously that $ \mathcal{A}$ is a subalgebra of $C_0(\mathbb{R}^{n})$ and $ \mathcal{A}$ separates points.
 Since $\mathbb{R}^{n}$ is a local compact Hausdorff space, $C_0(\mathbb{R}^{n})$ is a Banach algebra with the supremum norm  $
 \|f\|=\sup\{f(x): x\in\mathbb{R}^{n} \}$.  The Stone-Weierstrass theorem assures that $\mathcal{A}$ is dense in
 $C_0(\mathbb{R}^{n})$.

  Suppose that $0<p<\infty,$  and  $f\in L^p(\mathbb{R}^{n}).$ Let
 $$f_{N}(x)= f(x)\mathcal{X}_{\{x:|x|\leq N; |f(x)|\leq N\}}(x).$$
 It is clear to see that, $|f_{N}(x)|\leq N$, and
  $${\rm supp}f_N\subset \overline{B(0,N)}=\{x:\ |x|\leq N\}.$$
 Moreover, $|f_N-f|\leq 2|f|\in L^p(\mathbb{R}^{n})$, and, $\lim \limits _{N\rightarrow \infty}|f_N-f|^p = 0. $
 Thus, by the Lebesgue dominated convergence  theorem, we have
 $$
 \lim\limits _{N\rightarrow\infty}\int_{\mathbb{R}^{n}}|f_N-f|^{p}dx=\int_{\mathbb{R}^{n}}\lim\limits
 _{N\rightarrow\infty}|f_N-f|^{p}dx=0.
 $$
 Therefore, for $\varepsilon_0>0$, there exists an integer $N>1$, such that
\begin{equation}\label{13}
||f_N-f||_p^{p}<\varepsilon_0/3 .
\end{equation}
 Because $f_{N}(x)$ is a measurable function, according Lusin Theorem,  there exists a function $g_0\in C_0(\mathbb{R}^{n}),$
 such that
 $${\rm supp} g_0\subset B(0,N),\ \ \ |g_0(x)|\leq N,\ \ \ m(E_{N})< \varepsilon_0/3(2N)^p,$$
  where $E_{N}=\{ x:\ g_0(x)\neq f_{N}(x)\}$.
 Therefore,
 \begin{equation}\label{14}
 \begin{split}
\int_{\mathbb{R}^{n}}|g_0(x)-f_N|^{p}dx &= \int_{B(0,N)}|g_0(x)-f_N|^{p}dx \\
& \leq  (2N)^{p}\int_{E_{N}}dx< \varepsilon_0/3.
\end{split}
\end{equation}
Thus we obtain
\begin{equation}\label{15}
\|g_0-f\|_{p}^{p}= \|g_0-f_N+f_N-f\|_{p}^{p}\leq
\|g_0-f_N\|_{p}^{p}+\|f_N-f\|_{p}^{p}<\varepsilon_0/3+\varepsilon_0/3=2\varepsilon_0/3.
\end{equation}
Taking integers $\widetilde{l}_k$ such that $p\widetilde{l}_k>1,\ k=1,2,...,n,$
the fact ${\rm supp} g_0\subset B(0,N)$ implies
$$g_{0}(x)\prod \limits_{k=1}^{n}(1+x_{k}^{2})^{\widetilde{l}_k}\in C_0(\mathbb{R}^{n}).$$
Since $ \mathcal{A}$ is dense in $C_0(\mathbb{R}^{n})$, there exists
rational functions $r(x)\in \mathcal{A}$, such that
\begin{equation}\label{16}
\left|r(x)-g_{0}(x)\prod \limits_{k=1}^{n}(1+x_{k}^{2})^{\widetilde{l}_k}\right|<(\varepsilon_0/3)^{1/p}.
\end{equation}
 In fact, $r(x)$ can be written as
$$r(x)=\frac{P(x)}{\prod \limits_{k=1}^{n}(1+x_{k}^{2})^{l_k}},$$
where $P(x)=\sum_{(k)}\alpha_k x_{1}^{k_1}
 \cdot\cdot\cdot x_{n}^{k_n}, \ $\rm{deg}$_jP <2l_j, \ j=1,2,...,n$.
 Then, by \eqref{16},
$$\left|P(x)/\prod \limits_{k=1}^{n}(1+x_{k}^{2})^{l_k}-g_{0}(x)\prod
\limits_{k=1}^{n}(1+x_{k}^{2})^{\widetilde{l}_k}\right|<\varepsilon_0/3\pi^n,$$
or
$$\left|P(x)/\prod
\limits_{k=1}^{n}(1+x_{k}^{2})^{l_k+\widetilde{l}_k}-g_{0}(x)\right|<\frac{(\varepsilon_0/3)^{1/p}}{\pi^n\prod
\limits_{k=1}^{n}(1+x_{k}^{2})^{\widetilde{l}_k}}.$$
Obviously, $Q(x)=P(x)/\prod\limits_{k=1}^{n}(1+x_{k}^{2})^{l_k+\widetilde{l}_k}\in\mathcal{A}$. Thus,
$$\int_{\mathbb{R}^{n}}|Q(x)-g_{0}|^{p}dx
<\frac{1}{\pi^n}\int_{\mathbb{R}^{n}}\frac{(\varepsilon_0/3)}{\prod \limits_{k=1}^{n}(1+x_{k}^{2})^{p\widetilde{l}_k}}\ dx
<\varepsilon_0/3.$$
Therefore, for any $\varepsilon_0>0$, there exists $Q(x)\in\mathcal{A}$ such that
$$\|Q-f\|_{p}^{p}=\|Q-g_{0}+g_{0}-f\|_{p}^{p}\leq\|Q-g_{0}\|_{p}^{p}+\|g_{0}-f\|_{p}^{p}<\varepsilon_0.$$
Thus, for any $\varepsilon>0,$ taking $\varepsilon_k= \frac{\|f\|_{p}^{p}\varepsilon_0}{4^{k+3}},\ k=1,2,...,$ there exist
$Q_k(x)\in\mathcal{A}$, such that
$$\|Q_k-f\|_{p}^{p}< \varepsilon_k.$$
Moreover, the function $Q_k(z)$ is a rational function satisfying
$$\|Q_k\|_{p}^{p}=\|Q_k-f+f\|_{p}^{p}\leq\|Q_k-f\|_{p}^{p}+\|f\|_{p}^{p}<\varepsilon_k+\|f\|_{p}^{p}
=(1+\frac{\varepsilon}{4^{k+3}})\|f\|_{p}^{p}.$$
Thus the sequence of rational functions $Q_k(z)$ can be chosen such that
$$\|Q_k-Q_{k-1}\|_{p}^{p}=\|Q_k-f+f-Q_{k-1}\|_{p}^{p}\leq\|Q_k-f\|_{p}^{p}+\|f-Q_{k-1}\|_{p}^{p}<2\varepsilon_k,\ \ \
(k=2,3,\cdot\cdot\cdot).$$
Let
$$R_1(z)=Q_1(z),\  R_k(z)=Q_k(z)-Q_{k-1}(z),\ \ \ (k=2,3,\cdot\cdot\cdot).$$
Then $\{R_k(z)\}$ is a sequence of rational functions satisfying \eqref{6} and \eqref{7}. This completes the proof of  Lemma
3.2.

\begin{lemma} \ \ Suppose that   $0<p<1$ and $ R\in L^p(\mathbb{R}^{n})\bigcap \mathcal{A}$ is a rational function, where
$\mathcal{A}$ is the same as the $\mathcal{A}$ in Lemma 3.2.
 Then, there exist $2^n$ rational functions $R_{\sigma_j}(z)\in H^{p}(T_{\Gamma_{\sigma_j}}), \ j=1,2,...,2^n$, such that
 $$R(z)=\sum\limits _{j=1}^{2^n}R_{\sigma_j}(z),$$
  and
 \begin{equation}\label{8}
\sum\limits _{j=1}^{2^n}\|R_{\sigma_j}\|^p_{p}\leq C_{np}\|R\|^p_p,
 \end{equation}
where $C_{np}=2^n\left(\frac{2^{1-p}\pi}{1-p}\right)^{n}.$
\end{lemma}

{\bf Proof}\ \ \ Let $ R(x)\in L^p(\mathbb{R}^{n})\bigcap  \mathcal{A}$ be
a rational function.
 $R(z)$ can be written as
 $$
 R(z)=\frac{P(z)}{\prod\limits_{j=1}^n(1+z_j^{2})^{l_j}}=\frac{\sum_{(k)}\alpha_{k} z_{1}^{k_1}
 \cdot\cdot\cdot z_{n}^{k_n}}{\prod\limits_{j=1}^n(1+z_j^{2})^{l_j}},
 $$
 where $p(2l_j-{\rm deg}_jP)>1$ for all $j=1,2,...,n.$

 Let
 $$\beta(\xi)=\frac{i-\xi}{i+\xi}, \ \ \ \xi\in\mathbb{C}.$$
It is easy to know that $\beta(x)=e^{i \theta(x)},$ where $\theta(x)=\arg(i-x)-\arg(i+x)\in(-\pi,\pi)$ for $x\in\mathbb{R}$.

For each $\varphi=(\varphi_1,\varphi_2,...,\varphi_n)\in\mathbb{R}^{n}$, $z\in T_{\Gamma_{\sigma_j}}$.
  Rational functions $R_{\sigma_j}(z,\varphi) \ ( j=1,2,...,2^n)$ are defined as follows.

 $$R_{\sigma_j}(z,\varphi)=\frac{\prod\limits_{k=1}^{n}(\beta(z_k))^{m_k\sigma_j(k)}}
 {\prod\limits_{k=1}^{n}((\beta(z_k))^{m_k\sigma_j(k)}-e^{i\varphi_k\sigma_j(k)})}R(z),$$
 where $m_k>l_k+n\ (k=1,2,...,n)$ are positive integers .

 Then $R_{\sigma_j}(z,\varphi)$ can be written as,
  \begin{equation*}
 \begin{split}
R_{\sigma_j}(z,\varphi)& = \prod\limits_{\sigma_j(k)=+1}\frac{(\beta(z_k))^{m_k}}
 {(\beta(z_k))^{m_k}-e^{i\varphi_k}}\prod\limits_{\sigma_j(k)=-1}\frac{(\beta(z_k))^{-m_k}}
 {(\beta(z_k))^{-m_k}-e^{-i\varphi_k}}R(z)\\
&  = \prod\limits_{\sigma_j(k)=+1}\frac{(\beta(z_k))^{m_k}}
 {(\beta(z_k))^{m_k}-e^{i\varphi_k}}\prod\limits_{\sigma_j(k)=-1}\frac{e^{i\varphi_k}(\beta(z_k))^{m_k}}
 {(e^{i\varphi_k}-(\beta(z_k))^{m_k})(\beta(z_k))^{m_k}}R(z)\\
&  = \prod\limits_{\sigma_j(k)=+}\frac{(\beta(z_k))^{m_k}}
 {(\beta(z_k))^{m_k}-e^{i\varphi_k}}\prod\limits_{\sigma_j(k)=-1}\frac{(\beta(z_k))^{m_k}}
 {(\beta(z_k))^{m_k}-e^{i\varphi_k}}\prod\limits_{\sigma_j(k)=-1}\frac{-e^{i\varphi_k}}{(\beta(z_k))^{m_k}}R(z)\\
&  = \prod\limits_{k=1}^{n}\frac{(\beta(z_k))^{m_k}}
 {(\beta(z_k))^{m_k}-e^{i\varphi_k}}\prod\limits_{\sigma_j(k)=-1}\frac {-e^{i\varphi_k}}{(\beta(z_k))^{m_k}}R(z)\\
& = \frac{\prod\limits_{\sigma_j(k)=-1}\frac{-e^{i\varphi_k}}{(\beta(z_k))^{m_k}}}
  {\prod\limits_{k=1}^{n}\left (1-\frac{e^{i\varphi_k}}{(\beta(z_k))^{m_k}}\right )}R(z).\\
\end{split}
\end{equation*}
Therefore,
 $$\sum\limits _{j=1}^{2^n}R_{\sigma_j}(z,\varphi)=\frac{\sum\limits
 _{j=1}^{2^n}\prod\limits_{\sigma_j(k)=-1}\frac{-e^{i\varphi_k}}{(\beta(z_k))^{m_k}}}
  {\prod\limits_{k=1}^{n}\left (1-\frac{e^{i\varphi_k}}{(\beta(z_k))^{m_k}}\right )}R(z)= R(z).$$
Next we are to prove that $R_{\sigma_j}(z,\varphi)\in H^{p}(T_{\Gamma_{\sigma_j}}),$ for all $j=1,2,...,2^n.$

 Now we only consider the case that the base is the first octant of $\mathbb{R}^{n}$, because proofs of
 the other octants are similar.  For any $z\in T_{\Gamma_{\sigma_1}}$,
  \begin{equation}\label{18}
 \begin{split}
R_{\sigma_1}(z,\varphi)& = \frac{\prod\limits_{k=1}^{n}(\beta(z_k))^{m_k}}
 {\prod\limits_{k=1}^{n}((\beta(z_k))^{m_k}-e^{i\varphi_k})}R(z)\\
 & = \frac{\prod\limits_{k=1}^{n}(\frac{i-z_k}{i+z_k})^{m_k}}
 {\prod\limits_{k=1}^{n}((\beta(z_k))^{m_k}-e^{i\varphi_k})}\frac{P(z)}{\prod\limits_{k=1}^{n}(z_k-i)^{l_k}(z_k+i)^{l_k}}\\
&   = P(z)
 \prod\limits_{j=1}^{n}\frac{(-1)^{m_j}(z_j-i)^{m_j-l_j}}
 {(z_j+i)^{l_j+m_j}((\beta(z_j))^{m_j}-e^{i\varphi_j})}.
\end{split}
\end{equation}
Since $m_k > l_k,$ and $|\beta(z_k)|<1,\ |e^{i\varphi_k}|=1,$ for all $ k=1,2,...,n$, the function $R_{\sigma_1}(z)$ is a
rational function which is holomorphic in the tube $T_{\Gamma_{\sigma_1}}$.

Moreover, set
$$I_{\sigma_1} = \int_{\mathbb{R}^{n}}\int_{(-\pi, \pi)^{n}}|R_{\sigma_1}(x,\varphi)|^{p}\ d\varphi dx.$$
Then,
 \begin{equation}\label{25}
 \begin{split}
I_{\sigma_1}
& = \int_{\mathbb{R}^{n}}\int_{(-\pi, \pi)^{n}}\left|\frac{\prod\limits_{k=1}^{n}(\beta(x_k))^{m_k\sigma_1(k)}}
 {\prod\limits_{k=1}^{n}((\beta(x_k))^{m_k\sigma_1(k)}-e^{i\varphi_k\sigma_1(k)})}R(x)\right|^{p}\ d\varphi_1\cdot\cdot\cdot
 d\varphi_n dx\\
& = \int_{\mathbb{R}^{n}}\int_{-\pi}^{\pi}\cdot\cdot\cdot\int_{-\pi}^{\pi}\prod\limits_{k=1}^{n}\frac{|R(x)|^{p}}
 {|1-e^{i\varphi_k\sigma_1(k)-im_k\theta(x_k)}|^{p}}\ d\varphi_1\cdot\cdot\cdot d\varphi_n dx.\\
  \end{split}
\end{equation}
Observe that
\begin{equation}\label{226}
 \begin{split}
 &\int_{-\pi}^{\pi}\frac{1}
 {|1-e^{i\varphi_k\sigma_j(k)-im_k\theta(w_k)}|^{p}}\ d\varphi_k\\
 &=\int_{-\pi}^{\pi}\frac{1}
 {|1-e^{i\theta\sigma_j(k)}|^{p}}\ d\theta\\
 &=\int_{-\pi}^{\pi}\frac{d\theta} {2^{p}|\sin^{p}(\frac{\theta}{2})|}\leq \frac{4}{2^{p}}\int_{0}^{\frac{\pi}{2}}
 \frac{d\theta}
 {(\frac{2\theta}{\pi})^{p}}\leq \frac{2^{1-p}\pi}{1-p}.
  \end{split}
\end{equation}
By \eqref{25} and \eqref{226},
$$I_{\sigma_1}\leq \left(\frac{2^{1-p}\pi}{1-p}\right)^{n}\int_{\mathbb{R}^{n}}|R(x)|^{p}\ dx=
\left(\frac{2^{1-p}\pi}{1-p}\right)^{n}\|R\|_{p}^{p}.$$
Similarly,
$$I_{\sigma_k}\leq \left(\frac{2^{1-p}\pi}{1-p}\right)^{n}\|R\|_{p}^{p}$$
for $k=2,3,...,2^n.$

Therefore,
$$\sum\limits_{k=1}^{2^n}I_{\sigma_k}
=\int_{(-\pi, \pi)^{n}}\sum\limits_{k=1}^{2^n}\int_{\mathbb{R}^{n}}|R_{\sigma_k}(x,\varphi)|^{p}\ dxd\varphi\leq
2^n\left(\frac{2^{1-p}\pi}{1-p}\right)^{n}\|R\|^p_p.$$
 Thus, there exists a $\varphi=(\varphi_1,...,\varphi_n)\in (-\pi, \pi)^n$, such that
$$
\sum\limits_{k=1}^{2^n}\|R_{\sigma_k}\|^p_{p}\leq
2^n\left(\frac{2^{1-p}\pi}{1-p}\right)^{n}\|R\|^p_p,
$$
where $R_{\sigma_k}(x)=R_{\sigma_k}(x,\varphi),\ k=1,2,...,2^n$. This shows that the inequality \eqref{8} holds.
It is easy to know that
$$
\|R_{\sigma_k}\|^p_{p}\leq
2^n\left(\frac{2^{1-p}\pi}{1-p}\right)^{n}\|R\|^p_p,
$$
for all $ k=1,2,...,2^n$. So rational functions $R_{\sigma_k}(x)\in L^p(\mathbb{R}^n)$.  By Lemma 3.1, and that
$R_{\sigma_k}(z)$ is holomorphic in $T_{\Gamma_{\sigma_k}}$, we have
$$R_{\sigma_k}(z)\in H^p(T_{\Gamma_{\sigma_k}})$$
for all $ k=1,2,...,2^n$.

Thus, the proof of Lemma 3.3 is complete.

We still need the following lemmas.
\begin{lemma} (\cite{SW})\ \ Let $ B$ be an open cone in $\mathbb{R}^n$.
Suppose $F \in H^{p}(T_B),\ p > 0,$ and $B_0 \subset B$ satisfies $d(B_0, B^{c})= \inf \{|y_1-y_2|;\ y_1\in B_0,\ y_2 \in B^{c}\}\geq\varepsilon >0$, then there exists a constant $C_p(\varepsilon),$ depending on $\varepsilon$ and $p$  but not on $F$, such that
$$\sup \limits _{z\in T_{B_0}}|F(z)| \leq C_p(\varepsilon)\|F\|_{H^p}. $$
\end{lemma}
Given below offers a more precise estimation than that obtained in above Lemma 3.4.
\begin{lemma}\  Suppose that $f \in H^{p}(T_\Gamma),\ p > 0,$ and $ T_\Gamma$ is the tube with its base $\Gamma$ as the first octant of $\mathbb{R}^{n}$. If let
$f_\delta(z)= f(z+i\delta)$, for any $z=x+iy\in T_\Gamma$ and $\delta= (\delta_1,\delta_2,...,\delta_n)\in \Gamma$, then there holds
$$\sup \limits _{z\in T_{\Gamma}}|f_\delta(z)|\leq C_p \|f\|_{H^{p}}(\delta_1\cdot\cdot\cdot\delta_n)^{-\frac{1}{p}},$$
where $C_p= (\frac{2}{\pi})^{\frac{n}{p}}.$
\end{lemma}
We note that the proof of Lemma 3.5 is obtained in our recently work \cite{LDQ2}.
\subsection{ Proof of Theorem 3.1 }

Based on the above lemmas, we are now to prove Theorem 3.1.

{\bf Proof of Theorem 3.1}\ \ When $0<p<1$, by Lemma 3.2, for any $f(x)\in L^{p}(\mathbb{R}^{n})$, and  $\varepsilon>0,$ there exists
a sequence of rational functions $\{R_{k}(x)\}$, such that
 \begin{equation}\label{9}
  \sum_{k=1}^{\infty }||R_k||^p_p\leq (1+\varepsilon )\|f\|_p^p,
 \end{equation}
and
 \begin{equation}\label{10}
   \lim_{m\rightarrow \infty}||f-\sum_{k=1}^{m}R_k||_p
  =0.
 \end{equation}
  For each $k=1,2,...$, by Lemma 3.3, there exist $2^n$ rational functions $R_{k\sigma_j}(z)\in H^{p}(T_{\Gamma_{\sigma_j}})\ (
  j=1,2,...,2^n),$ such that
  \begin{equation}\label{11}
  R_k(z)=\sum\limits _{j=1}^{2^n}R_{k\sigma_j}(z),
  \end{equation}
  and
 \begin{equation}\label{12}
\sum\limits _{j=1}^{2^n}\|R_{k\sigma_j}\|^p_{p}\leq C_{np}\|R_k\|^p_p.
 \end{equation}
Therefore,
$$
  \sum_{k=1}^{\infty}\sum\limits_{j=1}^{2^{n}}||R_{k\sigma_j}||^p_p\leq \sum_{k=1}^{\infty}C_{np}\|R_k\|^p_p \leq
  (1+\varepsilon)C_{np}\|f\|_p^p,
$$
$$
   \lim_{m\rightarrow \infty}||f-\sum_{k=1}^{m}\sum\limits_{j=1}^{2^{n}}R_{k\sigma_j}||_p
  =0.
$$
By Lemma 3.1-2, there are
$$\|R_{k\sigma_j}\|^p_{H^p_{\sigma_j}}=\|R_{k\sigma_j}\|^p_{p},\ \ j= 1,2,...,2^{n},$$
which imply that the properties \eqref{2} and \eqref{3} hold.

Moreover, for any  $\delta=(\delta_1,\delta_2,...,\delta_n)\in \Gamma _{\sigma_j}$, and $ j=1,2,...,2^n$, by Lemma 3.4 and 3.5, there
exists a positive constant $ M_\delta< \infty $ such that
$$\left|R_{k\sigma_j}(z+i\delta)\right|\leq M_\delta\|R_{k\sigma_j}\|_{p}.$$
Hence,
$$\left|\sum_{k=1}^{m}R_{k\sigma_j}(z+i\delta)\right|^p\leq \sum_{k=1}^{m}|R_{k\sigma_j}(z+i\delta)|^p\leq
M_\delta\sum_{k=1}^{m}\|R_{\sigma_j}\|^p_{p},$$
for $z\in T_{\Gamma_{\sigma_j}}.$
This implies that the series $\sum\limits_{k=1}^{\infty}R_{k\sigma_j}(z)$ uniformly converges to a function $f_{\sigma_j}(z)$
in the tube domain
$T_{\Gamma_{\sigma_j}^\delta}=\{z=x+iy\in \mathbb{C}^n:\ x \in \mathbb{R}^n,\ \sigma_j(l)y_l>\sigma_j(l)\delta_l,\
l=1,2,...,n\}$ for any $\delta=(\delta_1,\delta_2,...,\delta_n)\in \Gamma _{\sigma_j}$.\\
 As a consequence, the function $f_{\sigma_j}(z)$ is holomorphic in $T_{\Gamma _{\sigma_j}}$.
Property \eqref{2} implies that property \eqref{4} holds.
By  properties of Hardy spaces on tubes, the non-tangential boundary limit $f_{\sigma_j}(x)$ of function $f_{\sigma_j}(z)\in H^{p}(T_{\Gamma_{\sigma_j}})$
exists for every $j=1,2,...,2^{n}$.\\
Therefore, property \eqref{3} implies that
$$f(x)=\sum\limits_{j=1}^{2^{n}}f_{\sigma_j}(x)$$
holds almost everywhere, and moreover,
$$\|f\|^{p}_{p} \leq \sum\limits_{j=1}^{2^{n}}\|f_{\sigma_j}\|^{p}_{p}\leq  A_p\|f\|_p^p.$$
Thus, the proof of Theorem 3.1 is complete.

\section{ Non-uniqueness of Hardy Space Decomposition}

 In this section, we are to answer the questions asked in the Remark 3.1 of Theorem 3.1. For the one dimension,
A.B. Aleksandrov (\cite{AB1}, \cite{CJA}) obtained the following theorem.

{\bf  Theorem B }\ (\cite{AB1} and \cite{CJA}) Let $0<p<1$ and $X^p$ denote the $L^p$ closure of the set of $f\in
L^p(\mathbb{R})$ which can be written in the form
$$f(x)= \sum\limits_{j=1}^N \frac{c_j}{x-a_j},\ \ a_j\in \mathbb{R}, \ c_j\in \mathbb{C}.$$
Then
$$X^p = H_+^p(\mathbb{R})\cap H_-^p(\mathbb{R}).$$

We note that, A.B. Aleksandrov's proof of Theorem M (\cite{AB1}, \cite{CJA}) is rather long involving vanishing moments and
the Hilbert transformation.  Deng and Qian  \cite{DQ} present a  more straightforward proof for Theorem M. In this section, our
aim is to extend Theorem M to higher dimensions. In order to do this, we need first to extend the following Theorem C obtained
by J.B. Garnett (\cite{Gar1}) and Theorem D obtained by Deng-Qian (\cite{DQ}) to higher dimensions.

{\bf  Theorem C } (\cite{Gar1}) Let $N$ be a positive integer. For $0<p<\infty,\ Np>1,$ the class $\mathfrak{w}_N$ is dense in
$H^p(\mathbb{C}^+)$, where $\mathfrak{w}_N$ is the family of $H^p(\mathbb{C}^+)$ functions satisfying

(i) $f(z)$ is infinitely differentiable in $\overline{\mathbb{C}^+}$,

(ii) $|z|^Nf(z)\rightarrow 0$ as $z \rightarrow \infty$.

{\bf  Theorem D  }(\cite{DQ}) \ Let $N$ be a positive integer. For $0<p<\infty,\ Np>1$, the class $\mathfrak{R}_N(i)$ is dense
in $H^p(\mathbb{C}^+)$, and the class $\mathfrak{R}_N(-i)$ is dense in $H^p(\mathbb{C}^-)$.
Where $\alpha\in \mathbb{C}$ and $\mathfrak{R}_N(\alpha)$ is the family of rational functions
$f(z)=(z+\alpha)^{-N-1}P(\frac{1}{z+\alpha}),$ $P(w)$ are polynomials.

We obtain the following three theorems for higher dimensions.
\begin{theorem}\label{Th5.3.1}
Let $N$ be a positive integer. For $0<p<\infty,\ pN>1$, the class $\mathfrak{w}_N$ is dense in $H^p(T_\Gamma),$ where
$\mathfrak{w}_N$ is the family of $H^p(T_\Gamma)$ functions satisfying

(i) $f(z)$ is infinitely differentiable in $\overline{T}_\Gamma$,

(ii) $|z|^Nf(z)\rightarrow 0$ as $|z|\rightarrow \infty$, where $|z|\rightarrow \infty$ means that $z_j\rightarrow \infty,\
1\leq j\leq n,\ z=(z_1,z_2,...,z_n)\in \overline{T}_\Gamma.$
\end{theorem}

{\bf Proof} We can approximate $f(z)\in H^p(T_\Gamma)$ by the smooth function $f(z+i/m)= f(z_1+i/m, ..., z_n+i/m)$.
In fact, the property that the existence  of the boundary limits functions of Hardy space functions assures that
$$\|f_m-f\|_{H^p}\rightarrow 0, \ \ m\rightarrow \infty.$$
We will construct the special functions $g_k(z)$ such that

(a) $g_k(z)\in \mathfrak{w}_N$,

(b) $|g_k(z)|\leq 1,\ \ z\in T_\Gamma$,

(c) $|g_k(z)|\rightarrow 1,\ \ z\in T_\Gamma,$ as $k\rightarrow \infty$.\\
Before we construct the above functions $g_k(z)$, we note that the functions
$$f_m(z)= g_m(z)f(z+i/m)$$
in $ \mathfrak{w}_N$ and then obtain the desired approximation.

 That is to say, if there exist the special functions $g_k(z)$, we can complete the proof of Theorem 4.1.

 As the heart of the proof, we are to construct the functions $g_k(z)$ in the following.
 Let $(w_1,w_2,...,w_n)\in T_{\Gamma},\ (\alpha_k,\alpha_k,...,\alpha_k)\in\mathbb{R}^n,\ 0<\alpha_k<1,$ and
 $\alpha_k\rightarrow 1$ as $k\rightarrow \infty.$\\
 Consider the function
 $$h_k(w)= \prod\limits_{j=1}^n\varphi_j(w_j),$$
 where $\varphi_j(w_j)=\left(\frac{w_j+\alpha_k}{1+\alpha_kw_j}\right)^{N+1},\ j=1,2,...,n$ has infinite $(N+1)-$fold zero at
 $-\alpha_k$.\\
Then,
 \begin{equation}\label{5.21}
|h_k(w)|= \left|\prod\limits_{j=1}^n\varphi_j(w_j)\right|
=\prod\limits_{j=1}^n\left|\frac{w_j+\alpha_k}{1+\alpha_kw_j}\right|^{N+1} < 1.
\end{equation}
Fixing $N$, $h_k(w)$ converges to $1$ uniformly on the compact set $\overline{D}\setminus \bigcup\limits_{k=1}^nE_k$, where
 $$D=D_1\times D_2\times\cdot\cdot\cdot\times D_n, \ \ D_k=|z|<1,\ k=1,2,...,n,$$
 and $$E_k=\{(w_1,...,w_{j-1},-\alpha_k,w_j,...,w_n):\ w_j\in D_k,\ k\neq j\}.$$
Then, for $w=(w_1,w_2,...,w_n)= \left(\frac{i-z_1}{i+z_1},...,\frac{i-z_n}{i+z_n}\right)$ and
$\alpha=(\alpha_k,\alpha_k,...,\alpha_k),$ we define the functions
$$g_k(z)= h_k(\alpha_k w).$$
Below we verify that the function $g_k(z)$ satisfies the three conditions (a) (b) and (c).\\
In fact, for condition (a), there is
\begin{equation*}
 \begin{split}
g_k(z) & = h_k(\alpha_k w)= \prod\limits_{j=1}^n\left(\frac{\alpha_k w_j+\alpha_k}{1+\alpha^2_kw_j}\right)^{N+1}\\
& = \prod\limits_{j=1}^n\alpha_k^{N+1}\left(\frac{ w_j+1}{1+\alpha^2_kw_j}\right)^{N+1}= \prod\limits_{j=1}^n\alpha_k^{N+1}\left(\frac{ \frac{i-z_j}{i+z_j}+1}{1+\alpha^2_k\frac{i-z_j}{i+z_j}}\right)^{N+1}\\
&= \prod\limits_{j=1}^n\alpha_k^{N+1}\left(\frac{2i}{(1+\alpha^2_k)i+(1-\alpha^2_k)z_j}\right)^{N+1}.
\end{split}
\end{equation*}
It is clearly that $g_k(z)$ satisfies the first condition $(i)$ of the class $\mathfrak{w}_N$.
Moreover, there holds
$|z|^Ng_k(z)\rightarrow 0$ as $|z|\rightarrow \infty$. This implies that $g_k(z)$ satisfies the second condition $(ii)$ of the
class $\mathfrak{w}_N$. Therefore, $g_k(z)\in \mathfrak{w}_N$, which shows that $g_k(z)$ satisfies the condition $(a)$.

For condition (b), from \eqref{5.21}, we can get that
$$|g_k(w)|= |h_k(\alpha_k w)|< 1.$$

 For (c), it is clear that

$$g_k(w)= h_k(\alpha_kw) \rightarrow 1,\ \ k\rightarrow \infty,\ z\in \overline{T}_\Gamma.$$ Thus, the proof  is complete.

We shall notice that the condition $Np>1$ implies that the above class $\mathfrak{w}_N$ is contained in $H^p(T_\Gamma)$. Let
$\alpha=(\alpha_1,\alpha_2,...,\alpha_n)\in \mathbb{C}^n$ and let $\mathfrak{R}_N(\alpha)$ be the family of the rational functions
$$f(z)=(z_1+\alpha_1)^{-N-1}(z_2+\alpha_2)^{-N-1}\cdot\cdot\cdot(z_n+\alpha_n)^{-N-1}
P(\frac{1}{z_1+\alpha_1},...,\frac{1}{z_n+\alpha_n}),$$
 where $z=(z_1,z_2,...,z_n)\in \mathbb{C}^n$ and $P(w)$ are polynomials. We notice that the class $\mathfrak{R}_N(\alpha)$ is
 contained in the class $\mathfrak{w}_N$ for Im$\alpha_j>0,\ j=1,2,...,n.$ Thus, we obtain the following results.

 \begin{theorem}\label{th5.3.2}
 Let $N$ be a positive integer. For $0<p<\infty$ and $ Np>1$, the class $\mathfrak{R}_N(i,i,...,i)$ is dense in
 $H^p(T_\Gamma).$
\end{theorem}

 \begin{corollary}
 Let $N$ be a positive integer. For $0<p<\infty$ and $ Np>1$, the class $\mathfrak{R}_N(\sigma_j(1)i,...,\sigma_j(n)i)$ is
 dense in $H^p(T_{\Gamma_{\sigma_j}}),\ j=1,2,...,2^n.$
 \end{corollary}
The proof of Corollary 4.1 is similar to the proof of Theorem 4.2, so we only  prove Theorem 4.2.\\
{\bf Proof of Theorem 4.2}
If $f(z)\in H^p(T_\Gamma),\ Np>1$, then, for any $\varepsilon>0,$ by Theorem 4.1, there exists function $f_N$ in
$H^p(T_\Gamma)\bigcap C^\infty(\overline{T}_\Gamma)$ such that
$$\lim\limits_{|z|\rightarrow \infty, z\in T_\Gamma}|z|^{N+1}f(z)=0,$$
and
$$\|f_N-f\|_{H^p}<\varepsilon.$$
The fractional linear mapping
$$z_j=\alpha(w_j)= i\frac{1-w_j}{1+w_j},\ \ j=1,2,...,n,$$
is a conformal mapping from the $n$-tuple unit disc
$$\mathbb{D}=\mathbb{D}_1\times \mathbb{D}_2 \times \cdot\cdot\cdot \mathbb{D}_n= \left\{(w_1,w_2,...,w_n)\in \mathbb{C}^n:
|w_j|<1,\ j=1,2,...,n\right\}$$
 to the first octant of $\mathbb{C}^n$, $T_\Gamma=\{(z_1,z_2,...,z_n)\in \mathbb{C}^n: {\rm Im}z_j>0,\ j=1,2,...,n\}$. \\
 Its inverse mapping is
$$w_j= \beta (z_j)= \frac{i-z_j}{i+z_j},\ \ j=1,2,...,n.$$
 Let
$$h_N(w)=f_N(\alpha(w_1),...,\alpha(w_n))$$
and $h_N(w_1,...,w_j,-1,w_{j+1},...,w_n)=0,\ \ j=1,2,...,n.$\\
Then $h_N(w)$ is continues in the closed disc $\overline{\mathbb{D}}$, and
$$h_N(w)\left|i\frac{1-w_1}{1+w_1}\right|^{N+1}\cdot\cdot\cdot\left|i\frac{1-w_n}{1+w_n}\right|^{N+1}\rightarrow 0, $$
as $ w_j\rightarrow w_0\in F,\ |w-w_0|\rightarrow 0,$ for $\ j=1,2,...,n$ and $w\in \overline{\mathbb{D}}\backslash  F$, where
$F=\bigcup\limits_{j=1}^nF_j,\ \ F_j=\{(w_1,...,w_j,-1,w_{j+1},...,w_n):|w_j|\leq 1,\  j=1,2,...,n\}.$

 Therefore, for $w_0\in F$,
$$\frac{h_N(w)}{\prod\limits_{j=1}^n(1+w_j)^{N+1}}\rightarrow 0,$$
as $\ w_j\rightarrow w_0,\ |w_j|\leq 1, \ j=1,2,...,n.$\\
Let
\begin{eqnarray}
\widetilde{h}_N(w) =
\begin{cases}
\frac{h_N(w)}{\prod\limits_{j=1}^n(1+w_j)^{N+1}}      & w \in \overline{\mathbb{D}}\backslash F \\
0   & w\in F
\end{cases}
\end{eqnarray}
Then $\widetilde{h}_N(w)$ is holomorphic in the $n$-tuple $\mathbb{D}$ and continues in the closed set $\overline{\mathbb{D}}$.\\
By the Stone-Weierstrass Theorem, there exists a polynomial $P_N(w)= P_N(w_1,...,w_n)$ such that
$$\left|\widetilde{h}_N(w)-P_N(w+1)\right|< \varepsilon,\ |w_j|\leq 1, \ w_j\neq -1,\ j=1,2,...,n,$$
that is
$$\left|h_N(w)/\prod\limits_{j=1}^n(1+w_j)^{N+1}-P_N(w+1)\right|< \varepsilon,\ |w_j|\leq 1, \ w_j\neq -1,\
j=1,2,...,n.$$
Thus,
$$\left|f_N(\alpha(w_1),...,\alpha(w_n))-\prod\limits_{j=1}^n(1+w_j)^{N+1}P_N(w+1)\right|\leq
\varepsilon\prod\limits_{j=1}^n|1+w_j|^{N+1},$$
where $ |w_j|\leq 1, \ w_j\neq 1,\ j=1,2,...,n.$\\
Since $z_j=\alpha(w_j)= i\frac{1-w_j}{1+w_j}$ and $w_j= \beta (z_j)= \frac{i-z_j}{i+z_j}$ for all $j=1,2,...,n$, the last
inequality becomes
$$\left|f_N(z)-\prod\limits_{j=1}^n\left(\frac{2i}{i+z_j}\right)^{N+1}P_N\left(\frac{2i}{i+z_n},...,\frac{2i}{i+z_n}\right)\right|\
\leq\varepsilon\prod\limits_{j=1}^n\left|\frac{2i}{i+z_j}\right|^{N+1},$$
for $z\in T_\Gamma,\ {\rm Im}z_j>0, \ j=1,2,...,n.$
Therefore, we can obtain that
$$\int_{\mathbb{R}^n}|f_N(x+iy)-R(x+iy)|^p\ dx<
\varepsilon^p2^{(N+1)p}\prod\limits_{j=1}^n\int_{-\infty}^{+\infty}\left|(1+x^2_j)\right|^{-(N+1)p} dx_j<\infty,$$
where $z=(z_1,z_2,...,z_n)\in T_\Gamma$ and
$$R(x+iy)= \prod\limits_{j=1}^n\left(\frac{2i}{i+z_j}\right)^{N+1}P_N\left(\frac{2i}{i+z_n},...,\frac{2i}{i+z_n}\right) \in
\mathfrak{R}_N(i,i,...,i).$$
This concludes that the class $\mathfrak{R}_N(i,i,...,i)$ is dense in $H^p(T_\Gamma).$ Therefore, the proof of Theorem 4.2 is complete.

The following result shows that the Hardy space decomposition of $L^p(\mathbb{R}^n)$ for $0<p<1$ is not unique and
the intersection space $\bigcap\limits_{j=1}^{2^n} H_{\sigma_j}^p(\mathbb{R}^n)$ is a non-empty set.

\begin{theorem}\label{th5.3.2}  Let $0<p<1$ and let $X^p$ denote the $L^p$ closure of the set of $L^p(\mathbb{R}^n)$-functions of the form
$$f(x_1,x_2,...,x_n)= \frac{P(x_1,...,x_n)}{\prod\limits_{k=1}^n\prod\limits_{j=1}^m(x_{k}-a_{kj})},$$
where $P(x_1,...,x_n)= \sum\limits_{s_1=0}^{l_1}\cdot\cdot\cdot\sum\limits_{s_n=0}^{l_n}a_{s_1...s_n}x_1^{s_1}\cdot\cdot\cdot
x_n^{s_n},\ l= \max\{l_1,...,l_n\}, $ $(m-l)p>1,$ and  $a_{kj}\neq a_{km} \in \mathbb{R}$ for $ j\neq m$ and $ k=1,2,...,n.$ \\
Then
$$X^p= \bigcap\limits_{j=1}^{2^n} H_{\sigma_j}^p(\mathbb{R}^n)\neq \emptyset.$$
\end{theorem}
{\bf Proof}  Firstly, we will show that $X^p\subseteq\bigcap\limits_{j=1}^{2^n} H_{\sigma_j}^p(\mathbb{R}^n).$
 Let
 $$f(z_1,z_2,...,z_n)=\frac{\sum\limits_{s_1=0}^{l_1}\cdot\cdot\cdot\sum\limits_{s_n=0}^{l_n}a_{s_1...s_n}z_1^{s_1}\cdot\cdot\cdot
 z_n^{s_n}}{\prod\limits_{k=1}^n\prod\limits_{j=1}^m(z_{k}-a_{kj})},\ \ $$
  $a_k=(a_{k1},...,a_{km})\in \mathbb{R}^n,\ k=1,2,...,n.$\\
Then $f(z)$ is a rational function with singular part being contained  in $\bigcup\limits_{l=1}^{n}G_l$, where $G_l=\{(z_1,...,z_n):\
{\rm Im}z_l=0,\ z_k\in\mathbb{C},\ k\neq l\}$. It is clear that the set $G_j$ is not in the octant $T_{\Gamma_{\sigma_j}}$ and
$\bigcup\limits_{l=1}^{n}G_l\subseteq \bigcap\limits_{j=1}^{2^n}\partial T_{\Gamma_{\sigma_j}}$. So $f(z)$ is holomorphic in
$T_{\Gamma_{\sigma_j}}$, for all  $ j=1,2,...,2^n$.

Moreover,
\begin{equation}\label{6.3.2}
 \begin{split}
\int_{\mathbb{R}^n}|f(x+iy)|^p\ dx & =
\int_{\mathbb{R}^n}\frac{\left|\sum\limits_{s_1=0}^{l_1}\cdot\cdot\cdot\sum\limits_{s_n=0}^{l_n}a_{s_1...s_n}z_1^{s_1}
\cdot\cdot\cdot z_n^{s_n}\right|^p}{\left|\prod\limits_{k=1}^n\prod\limits_{j=1}^m(z_{k}-a_{kj})\right|^p}\ dx\\
& \leq \sum\limits_{s}|a_{s_1...s_n}|^p\prod\limits_{k=1}^n\int_{-\infty}^{+\infty}
\frac{\left|z_k^{s_k}\right|^p}{\left|\prod\limits_{j=1}^m(z_{k}-a_{kj})\right|^p}\ dx_k.
\end{split}
\end{equation}
Observe that,  the function
$$R_k(z_k)=\frac{z_k^{s_k}}{\prod\limits_{j=1}^m(z_{k}-a_{kj})},$$
is holomorphic in the upper-half and the lower-half complex planes.
 Moreover, we can prove that $R_k(x_k)\in L^p(\mathbb{R})$. In fact,
\begin{equation*}
 \begin{split}
\int_{-\infty}^{+\infty}|R_k(x_k)|^p\ dx_k & =
\int_{-\infty}^{+\infty}\frac{\left|x_k^{s_k}\right|^p}{\left|\prod\limits_{j=1}^m(x_{k}-a_{kj})\right|^p}\ dx_k \\
&=\int_{|x_k|\leq M}\frac{\left|x_k^{s_k}\right|^p}{\left|\prod\limits_{j=1}^m(x_{k}-a_{kj})\right|^p}\ dx_k +\int_{|x_k|>
M}\frac{\left|x_k^{s_k}\right|^p}{\left|\prod\limits_{j=1}^m(x_{k}-a_{kj})\right|^p}\ dx_k.
\end{split}
\end{equation*}
where $M$ is sufficiently large so that the interval $(-M,\ M)$ contains all the poles of $R_k(x_k)$.
Then, for the first integral in the right hand of the last inequality,
$$\int_{|x_k|\leq M}\frac{\left|x_k^{s_k}\right|^p}{\left|\prod\limits_{j=1}^m(x_{k}-a_{kj})\right|^p}\ dx_k \leq
M^{s_kp}\int_{|x_k|\leq M}\frac{ dx_k}{\prod\limits_{j=1}^m|x_{k}-a_{kj}|^p},$$
which is finite since $0<p<1.$\\
For the other integral, we have
\begin{equation*}
 \begin{split}
\int_{|x_k|> M}\frac{\left|x_k^{s_k}\right|^p}{\left|\prod\limits_{j=1}^m(x_{k}-a_{kj})\right|^p}\ dx_k &= \int_{|x_k|>
M}\frac{\left|x_k^{s_k}\right|^p}{\left|x^m_{k}\prod\limits_{j=1}^m(1-a_{kj}/x_{k})\right|^p}\ dx_k\\
& \leq \int_{|x_k|> M}\frac{dx_k}{\prod\limits_{j=1}^m|1-a_{kj}/x_{k}|^p|x|^{(m-s_k)p}_{k}}\\
& \leq \prod\limits_{j=1}^m\frac{\left|M-a_{kj}\right|^p}{\left|M\right|^p}\int_{|x_k|> M}\frac{dx_k}{|x|^{(m-s_k)p}_{k}},
\end{split}
\end{equation*}
which is also finite since $(m-s_k)p>1.$ Therefore, $R_k(x_k)\in L^p(\mathbb{R}).$ \\
By Lemma 3.1, together with $R_k(x_k)\in L^p(\mathbb{R})$ and the analyticity of $R_k(x_k)$, we get that
$R_k(z_k)\in H^p(\mathbb{C}^{\pm}).$
Thus, $\|R\|_{H^p}\leq \|R\|_{L^p};$ and the relation \eqref{6.3.2} becomes
$$\int_{\mathbb{R}^n}|f(x+iy)|^p\ dx \leq \sum\limits_{s}|a_{s_1...s_n}|^p\|R\|^{np}_{L^p} < \infty.$$
Hence, we have
$$f(z)\in \bigcap\limits_{j=1}^{2^n}H^p(T_{\Gamma_{\sigma_j}}),\ \ j=1,2,...,2^n.$$
This shows that $$f(x)\in \bigcap\limits_{j=1}^{2^n} H_{\sigma_j}^p(\mathbb{R}^n).$$
Hence $X^p\subseteq\bigcap\limits_{j=1}^{2^n} H_{\sigma_j}^p(\mathbb{R}^n)$ as desired.

Next, we will show   $X^p\supseteq\bigcap\limits_{j=1}^{2^n} H_{\sigma_j}^p(\mathbb{R}^n).$

Let there exist $f_{\sigma_j}(z)\in H^p(T_{\Gamma_{\sigma_j}}), \ j=1,2,...,2^n,$ such that
$$f(x)=f_{\sigma_j}(x)=f_{\sigma_l}(x),$$
 for all $1\leq j,\ l\leq2^n,$ a.e. $x\in \mathbb{R}^n.$\\
By Theorem 4.2 and Corollary 4.1, for any $\varepsilon>1$, there exist
$$R_{\sigma_j}\in \mathfrak{R}_N(\sigma_j i)$$
such that
$$\|f_{\sigma_j}-R_{\sigma_j}\|^p_{H^p_{\sigma_j}}= \|f-R_{\sigma_j}\|^p_{p}< \frac{\varepsilon}{4},$$
for $j=1,2,...,2^n.$ \\
The fact
$R_{\sigma_j}\in \mathfrak{R}_N(\sigma_j i)$ implies that there exist polynomials $P_{\sigma_j},\ j=1,2,...,2^n,$
such that for  $\beta (z_k)= \frac{i-z_k}{i+z_k},\ \ k=1,2,...,n,$
$$R_{\sigma_j}(z)=P_{\sigma_j}((\beta(z_1))^{\sigma_j(1)}+1,...,(\beta(z_n))^{\sigma_j(n)}+1)
\prod\limits_{k=1}^n\left(\beta^{\sigma_j(k)}(z_k)+1\right)^{N+1},$$
where $(N+1)p>1, \ N+1> \max\{{\rm deg}P_{\sigma_j}:\ j=1,2,...,2^n\}.$

It is easy to see that the singular part of $R_{\sigma_j}(z)$ are contained in $\bigcup\limits_{k=1}^nH_k$ which is not in
$T_{\Gamma_{\sigma_j}}$, where $H_k=\{(z_1,...,z_n):\ z_k=-\sigma_j(k)i, \ z_l\in\mathbb{C}, \ k\neq l\}$. So $R_{\sigma_j}(z)$
is a holomorphic rational function in $T_{\Gamma_{\sigma_j}},\ j=1,2,...,2^.$ \\
Let
$$R(z,\varphi)=R_{\sigma_1}(z)+ \frac{\sum\limits_{j=1}^{2^n}\prod\limits_{k=1}^n(-e^{i\varphi_k})^{\frac{1+\sigma_j(k)}{2}}
\prod\limits_{k=1}^n(\beta(z_k))^{\frac{1-\sigma_j(k)}{2}m}(R_{\sigma_j}(z)-R_{\sigma_1}(z))}
{\prod\limits_{k=1}^n((\beta(z_k))^m-e^{i\varphi_k})},$$
where  $(m-\max\{{\rm deg}P_{\sigma_j}:\ j=1,2,...,2^n\}-N-1)p>1$, and $\varphi=(\varphi_1,...,\varphi_n)\in \mathbb{R}^n$. We thus are aware that the singular part of the rational function $R(z,\varphi)$ is contained in $\bigcap\limits_{j=1}^{2^n}\partial
T_{\Gamma_{\sigma_j}}=\mathbb{R}^n$. So $R(z,\varphi)$ is holomorphic in $\bigcup\limits_{j=1}^{2^n}T_{\Gamma_{\sigma_j}}$.\\
Note that $\beta(x_k)=e^{i\theta(x_k)}$, where $\theta(x_k)= \arg(i_k-x_k)-\arg(i_k+x_k)\in [-\pi, \pi]$ for $x_k\in \mathbb{R}$.
Moreover, set
$$J = \int_{\mathbb{R}^n}\int_{[-\pi, \pi]^n}|R(x,\varphi)-R_{\sigma_1}(x)|^p\ d\varphi dx.$$
Then
\begin{equation*}
 \begin{split}
J & = \int_{\mathbb{R}^n}\int_{[-\pi, \pi]^n}
\frac{\left|\sum\limits_{j=1}^{2^n}\prod\limits_{k=1}^n(-e^{i\varphi_k})^{\frac{1+\sigma_j(k)}{2}}
\prod\limits_{k=1}^n(\beta(z_k))^{\frac{1-\sigma_j(k)}{2}m}(R_{\sigma_j}(x)-R_{\sigma_1}(x))\right|^p}
{\left|\prod\limits_{k=1}^n((\beta(z_k))^m-e^{i\varphi_k})\right|^p}\ d\varphi dx\\
& \leq \sum\limits_{j=1}^{2^n}\int_{\mathbb{R}^n}\int_{[-\pi,
\pi]^n}|R_{\sigma_j}(x)-R_{\sigma_1}(x)|^p\prod\limits_{k=1}^n\left|
\frac{(\beta(z_k))^{\frac{1-\sigma_j(k)}{2}m}}
{(\beta(z_k))^m-e^{i\varphi_k}}\right|^p\ d\varphi dx\\
& = \sum\limits_{j=1}^{2^n}\int_{\mathbb{R}^n}\int_{[-\pi, \pi]^n}
\frac{|R_{\sigma_j}(x)-R_{\sigma_1}(x)|^p d\varphi dx}
{\prod\limits_{\sigma_j(k)=1}|e^{im\theta(x_k)}-e^{i\varphi_k}|^p\prod\limits_{\sigma_j(k)=-1}
|1-e^{i\varphi_k-im\theta(x_k)}|^p}.
\end{split}
\end{equation*}
Observe that
\begin{equation}\label{26}
 \begin{split}
 &\int_{-\pi}^{\pi}\frac{1}
 {|1-e^{i\varphi_k-im\theta(w_k)}|^{p}}\ d\varphi_k\\
 &=\int_{-\pi}^{\pi}\frac{1}
 {|1-e^{i\theta}|^{p}}\ d\theta\\
 &=\int_{-\pi}^{\pi}\frac{d\theta}
 {2^{p}|\sin^{p}(\frac{\theta}{2})|}\leq \frac{4}{2^{p}}\int_{0}^{\frac{\pi}{2}}\frac{d\theta}
 {(\frac{2\theta}{\pi})^{p}}\leq \frac{2^{1-p}\pi}{1-p},
  \end{split}
\end{equation}
and similarly,
\begin{equation}\label{26}
 \begin{split}
 &\int_{-\pi}^{\pi}\frac{1}
{|e^{im\theta(x_k)}-e^{i\varphi_k}|^p}\ d\varphi_k\\
 &=\int_{-\pi}^{\pi}\frac{d\theta}
 {2^{p}|\sin^{p}(\frac{\theta}{2})|}\leq \frac{4}{2^{p}}\int_{0}^{\frac{\pi}{2}}\frac{d\theta}
 {(\frac{2\theta}{\pi})^{p}}\leq \frac{2^{1-p}\pi}{1-p}.
  \end{split}
\end{equation}
From the above relations, we obtain
$$J \leq \left(\frac{2^{1-p}\pi}{1-p}\right)^n\int_{\mathbb{R}^n}|R_{\sigma_j}(x)-R_{\sigma_1}(x)|^p\ dx.$$
Therefore, there is a real vector $\varphi\in \mathbb{R}^n$ such that
$$J \leq \left(\frac{2^{1-p}\pi}{1-p}\right)^n(\|R_{\sigma_j}-f\|^p_p+\|f-R_{\sigma_1}\|^p_p)
\leq \left(\frac{2\pi}{1-p}\right)^n\varepsilon.$$
Thus, we have
\begin{equation}\label{53}
 \begin{split}
&\int_{\mathbb{R}^n}|R(x,\varphi)-f(x)|^p\ dx \\
& \leq \int_{\mathbb{R}^n}|R(x,\varphi)-R_{\sigma_j}(x)|^p\ dx +\int_{\mathbb{R}^n}|R_{\sigma_j}(x)-f(x)|^p\ dx \\
& \leq \left(\frac{2\pi}{1-p}\right)^n\varepsilon +\frac{\varepsilon}{4}.
\end{split}
\end{equation}
Hence, $R(x,\varphi) \in L^p(\mathbb{R}^n)$. This, together with the definition of $R(x,\varphi)$, implies that
$R(x,\varphi)\in X^p.$
Therefore, $f(x)\in X^p$. The proof is complete.


\end{document}